\documentclass[12pt]{article}
\usepackage{amsmath, amsthm, amssymb}
\numberwithin{equation}{section}
\newtheorem{theorem}{Theorem}

\begin{document}
\author{Ajai Choudhry}
\title{Rational quadrilaterals}
\date{}
\maketitle

\begin{abstract} 
A quadrilateral is said to be rational if its four sides, the two diagonals and the area are all expressible by rational numbers. The problem of constructing rational quadrilaterals dates back to the seventh century when Brahmagupta gave an elegant solution of  the problem. In 1848 Kummer gave a method of generating all rational quadrilaterals. In this paper we present an alternative method of generating all rational quadrilaterals. For rational cyclic quadrilaterals,  we obtain  a complete parametrization and for  noncyclic rational quadrilaterals, we give several parametrizations in terms of quadratic and quartic polynomials. The parametrizations obtained in this paper are simpler than the known parametrizations of rational quadrilaterals.  We also describe how further parametrizations of rational quadrilaterals may be obtained.

\end{abstract}

\noindent Mathematics Subject Classification 2020: 11D09
\smallskip
  
\noindent Keywords: Rational quadrilateral; cyclic quadrilateral; Brahmagupta 

\hspace*{0.54in} quadrilateral.

\section{Introduction}
Since ancient times there has been considerable interest in geometric objects such as triangles, quadrilaterals and other polygons such that the lengths of their  sides and  diagonals, as well as their  areas are expressible by rational numbers. This paper is concerned with rational quadrilaterals --- a rational  quadrilateral being defined as one whose  four sides, the two diagonals and the area are  given by  rational numbers.

In the seventh century Brahmagupta gave an elegant method of constructing rational quadrilaterals. Brahmagupta's method was elaborated further by Bhaskara in the 12$^{\rm th}$ century and by Chasles in the 19$^{\rm th}$ century (as quoted by Dickson \cite[pp. 216--217]{Di2}). In 1848 Kummer \cite{Ku} gave a method of generating all rational quadrilaterals. The method, however, leads to a cumbersome final condition that requires a parametrized quartic function be made a perfect square. The complete solution of this final condition is not yet known. In 1921 Dickson \cite{Di1} presented Kummer's construction in a somewhat simplified manner.  

Various  related problems concerning rational quadrilaterals have been studied by several  mathematicians (\cite{Al1}, \cite{Al2}, \cite{BM}, \cite[pp. 216--221]{Di2}, \cite{HS}, \cite{IKM}, \cite{Mo1}, \cite{Sa}). Some authors have   adopted a relaxed definition of rational quadrilaterals that requires only the four  sides and the two diagonals of the quadrilateral to be rational  without the additional condition of the area being rational while some others have studied only cyclic quadrilaterals. For instance, Sastry \cite[pp.\ 170--171]{Sa} has given a parametrization of cyclic quadrilaterals whose sides, diagonals and the area are all expressible by rational numbers. However, regarding the original problem of constructing quadrilaterals with rational sides, diagonals and area, there has been no progress in over a hundred years since Dickson's  simplification in 1921 of the method described by Kummer.

In this paper we describe a method, different from that of Kummer,  of generating all  rational convex quadrilaterals. We first given in Section 2 a necessary and sufficient condition for the existence of rational quadrilaterals. In Section 3 we obtain several examples of rational quadrilaterals whose sides, diagonals and area are given by multivariate polynomials. These parametrizations   are more general and  simpler than the examples obtained by Kummer's method. In fact, as a special case, we obtain, for  all rational cyclic quadrilaterals, a complete parametrization that is simpler than the known parametrizations of such quadrilaterals. 

\section{A necessary and sufficient condition for rational quadrilaterals}\label{necsuffcond}
We now prove  a theorem that describes a necessary and sufficient condition for the existence of a rational convex quadrilateral.

\begin{theorem}\label{Thcond} A necessary and sufficient condition for the existence of a rational convex quadrilateral is that there exist  rational numbers $a, b, c, d, e, f, x_1$, $x_2, y_1$ and $y_2$ satisfying the simultaneous diophantine equations,
\begin{align}
x_1^2+y_1^2&=a^2,\label{cond1} \\ 
(e-x_1)^2+y_1^2&=b^2, \label{cond2} \\
(e-x_2)^2+y_2^2&=c^2, \label{cond3} \\
x_2^2+y_2^2&=d^2, \label{cond4} \\
(x_1-x_2)^2+(y_1-y_2)^2&=f^2. \label{cond5} 
\end{align}
and such that the following rational numbers are all positive: $a, b, c, d, e, f$, $-y_1y_2 , -(y_1 - y_2)(x_1y_2 - x_2y_1)$ and $(y_1 - y_2)(ey_1 - ey_2 + x_1y_2 - x_2y_1)$.
\end{theorem}

\begin{figure}[h]
\begin{center}
\setlength{\unitlength}{0.5cm}
\begin{picture}(9,9)(1,-5)
\thicklines
\put(0,0){\line(3,2){2.4}}
\put(2.4,1.6){\line(5, -2){4.0}}
\put(6.4, 0){\line(-1,-2){2.4}}
\put(0,0){\line(5,-6){4.0}}
\thinlines
\put(0,0){\line(1,0){9.4}}
\put(2.4,0){\line(0,1){1.6}}
\put(4,0){\line(0,-1){4.8}}
\put(0,0){\line(0,1){2.0}}
\put(0,0){\line(0,-1){5.5}}
\put(9.5,-0.2){{\it \footnotesize X}}
\put(-0.2, 2.2){{\it \footnotesize Y}}
\put(-0.7,-0.2){{\it \footnotesize O}}
\put(6.4, 0.2){{\it \footnotesize B}\,{\footnotesize($e, 0$)}}
\put(1.9, -0.6){{ \it \footnotesize D}}
\put(2.2, 1.7){{\it \footnotesize A}{\footnotesize($x_1, y_1$)}}
\put(3.7, -5.5){{\it \footnotesize C}\,{\footnotesize($x_2, y_2$)}}
\put(3.9, 0.1){{\it \footnotesize E}}
\put(0.68, 0.98){{ \it \footnotesize a}}
\put(4.3, 0.85){{ \it \footnotesize b}}
\put(4.8, -3.1){{ \it \footnotesize c}}
\put(1.4, -3.1){{ \it \footnotesize d}}
\end{picture}
\caption{Rational quadrilateral}
\end{center}
\end{figure}
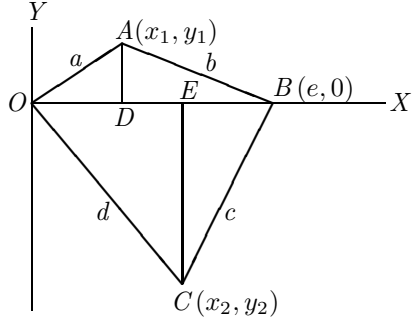
\begin{proof}
Without loss of generality, let $OABC$ be a rational convex quadrilateral with one   vertex $O$ at the origin of the $x$-$y$ plane and with its diagonal $OB$ along the $x$-axis as shown in Figure 1. Let the lengths of the sides $OA, AB, BC$ and $OC$ of the quadrilateral be $a, b, c$ and $d$ respectively, and let the lengths of the diagonals $OB$ and $AC$ be $e$ and $f$ respectively.  Let $D$ and $E$ be the feet of the perpendiculars to the diagonal $OB$ drawn from the vertices $A$ and $C$ respectively.   Let the coordinates of the vertices $A$ and $C$ be $(x_1, y_1)$ and $(x_2, y_2)$ respectively.  It readily follows from the theorem of Pythagoras that the relations \eqref{cond1}--\eqref{cond5} must be satisfied. 

Since  $OABC$ is assumed to be a rational quadrilateral, $a, b, c, d, e$ and $f$ are rational numbers and the area $A$ of the quadrilateral is rational. We will now prove that the four numbers $x_i, y_i, i=1, 2$, are all rational. Since all the three sides of the two triangles $OAB$  and $OCB$  are rational, the cosines of the angles $AOB$ and $COB$ are rational, and hence the lengths $OA \cos AOB$ and $OC \cos COB$, that is, $x_1$ and $x_2$  are rational. 

Further, it follows from Heron's formula for the area of a triangle that the areas of the two triangles $OAB$  and $OCB$ may be written as $\sqrt{r_1}$ and $\sqrt{r_2}$ where $r_1$ and $r_2$ are rational numbers. Since, by definition, the area $A$ of the quadrilateral  $OABC$ is rational, it follows that $\sqrt{r_1}+\sqrt{r_2}=A$, or,  $\sqrt{r_1} =A-\sqrt{r_2}$, or, $r_1=A+r_2-2A\sqrt{r_2}$, and hence $\sqrt{r_2}$ is a rational number. Similarly, $\sqrt{r_1}$ is also a rational number. Thus, the areas of the two triangles $OAB$  and $OCB$ are rational. Since we may write these areas as $OB \times AD/2$ and $OB \times CE/2$ respectively, and $OB$ is of rational length, it follows that the lengths of the line segments $AD$ and $CE$ are rational, that is, both $y_1$ and $y_2$ are rational numbers. 

 We will now prove that all the rational numbers mentioned towards the end of the theorem are positive. We first note that $a, b, c, d, e, f$, being the lengths of the sides and diagonals of the quadrilateral, are necessarily positive rational numbers. Further, since  $OABC$ is a convex quadrilateral, the diagonal $OB$ lies wholly within the quadrilateral, and hence the vertices $A$ and $C$ must be on opposite sides of the diagonal $OB$, and hence the ordinates $y_1$ and $y_2$ must be of opposite signs, and thus $-y_1y_2$  has to be positive. Finally, since $OABC$ is a convex quadrilateral, the point of intersection $M$ of the two diagonals must strictly lie somewhere between the two vertices $O$ and $B$, that is, the abscissa of the point $M$ must be between $0$ and $e$.  The coordinates of the point $M$ are  readily worked out to be $(-(x_1y_2 - x_2y_1)/(y_1 - y_2), 0)$, hence we must have
\[
0 < -(x_1y_2 - x_2y_1)/(y_1 - y_2) < e.
\]
It now readily follows that the numbers $-(y_1-y_2)(x_1y_2 - x_2y_1)$ and $(y_1-y_2)(ey_1 - ey_2 + x_1y_2 - x_2y_1)$ are both positive. 

We have thus proved that if there exists a rational convex quadrilateral, then there exist rational numbers $a, b, c, d, e, f, x_1, x_2, y_1$ and $y_2$ satisfying the relations \eqref{cond1}--\eqref{cond5} and such that the numbers $a, b, c, d, e, f, -y_1y_2 , -(y_1 - y_2)(x_1y_2 - x_2y_1)$ and $(y_1 - y_2)(ey_1 - ey_2 + x_1y_2 - x_2y_1)$ are all positive. We have thus shown that the condition stated in the theorem is necessary. 

Next let there exist rational numbers $a, b, c, d, e, f, x_1, x_2, y_1$ and $y_2$ satisfying the  conditions of the theorem. We will construct a quadrilateral $OABC$ with vertices $O, A, B$ and $C$ at the points $(0, 0), (x_1, y_1), (e, 0)$ and $(x_2, y_2)$ respectively as shown in Figure 1. In view of the relations \eqref{cond1}--\eqref{cond5}, it is clear that the lengths of the four sides of the quadrilateral are $a, b, c$, and $d$ respectively, and the two diagonals $OB$ and $AC$ have lengths $e$ and $f$ respectively, and hence these lengths are all rational. 

Further, the line segments $OB$, $AD$ and $CE$ have rational lengths, hence both the  triangles $OAB$ and $OCB$  have rational areas, and thus the area of the quadrilateral $OABC$ is rational. Finally it is readily seen that the conditions  that the numbers  $-y_1y_2 , -(y_1 - y_2)(x_1y_2 - x_2y_1)$ and $(y_1 - y_2)(ey_1 - ey_2 + x_1y_2 - x_2y_1)$ are all positive, ensure that the quadrilateral is convex.  Thus, the condition stated in the theorem is sufficient for the existence of a rational convex quadrilateral. This proves the theorem.
\end{proof}

It is readily seen that if 
\[
(a, b, c, d, e, f, x_1, x_2, y_1, y_2)=(a^{\prime}, b^{\prime}, c^{\prime}, d^{\prime}, e^{\prime}, f^{\prime}, x_1^{\prime}, x_2^{\prime}, y_1^{\prime}, y_2^{\prime})\]
is a solution of the simultaneous diophantine Eqs. \eqref{cond1}--\eqref{cond5}, then 
$(\varepsilon_1 a^{\prime},\varepsilon_2 b^{\prime}$, $\varepsilon_3 c^{\prime}, \varepsilon_4 d^{\prime}, \varepsilon_5 e^{\prime}, \varepsilon_6 f^{\prime}, \varepsilon_5 x_1^{\prime}, \varepsilon_5 x_2^{\prime}, \varepsilon_7 y_1^{\prime}, \varepsilon_7y_2^{\prime})$
is also a solution where the value of each $\varepsilon_i, i=1, \ldots, 7$, is taken  either as $+1$ or $-1$. Accordingly, if we obtain a numerical solution of Eqs.  \eqref{cond1}--\eqref{cond5} such that all the conditions of a rational convex quadrilateral are not satisfied, it is sometimes, though not always, possible to choose appropriate values of $\varepsilon_i$ and obtain a solution that satisfies all the conditions of a convex quadrilateral. 

\section{Parametrizations of  rational quadrilaterals}\label{paramsols}
We will now solve the simultaneous diophantine Eqs. \eqref{cond1}--\eqref{cond5}. It is not difficult to obtain the complete solution of any four of the five Eqs. \eqref{cond1}--\eqref{cond5}, and in each case the final equation reduces to a parameterized quartic function to be made a perfect square. While rational solutions of this final condition can be obtained, these are generally quite cumbersome. This situation is similar to the one that arises on following Kummer's method. 

In the next section we will solve  the four Eqs. \eqref{cond1}--\eqref{cond4} in a manner such that the final condition arising from Eq. \eqref{cond5} has simple solutions.

\subsection{Quadrilaterals with four sides, one diagonal and the area expressible by rational numbers}\label{solcond1234}
We will now obtain the complete solution of  the four simultaneous diophantine Eqs. \eqref{cond1}--\eqref{cond4}. We will repeatedly make use of the fact that the complete solution of the diophantine equation $X_1X_2=Y_1Y_2$ is given by $X_1=gm, X_2=hn, Y_1=gn, Y_2=hm$, where $g, h, m $ and $n$ are arbitrary parameters (see \cite[p.\ 69]{Si}). The complete solution of Eqs. \eqref{cond1}--\eqref{cond4}  will also yield a parametrization of all quadrilaterals $OABC$ in which the four sides, the diagonal $OB$ and the area of the quadrilateral $OABC$ are  rational. 

Taking the difference of the two equations, \eqref{cond1} and \eqref{cond2}, we get $-e(e-2x_1) = (a + b)(a - b)$, hence we may write  $-e=2u_1v_1, e-2x_1 =2u_2v_2,  a + b=2u_1v_2, a - b=2u_2v_1$, when we get the solution,
\begin{equation}
a = u_1v_2 + u_2v_1,\quad  b = u_1v_2 - u_2v_1, \quad e = -2u_1v_1,\quad  x_1 = -u_1v_1 - u_2v_2, \label{soleq12}
\end{equation}
where $u_1, u_2, v_1, v_2$ are arbitrary parameters.

Similarly, the complete solution of the difference of the two equations \eqref{cond3} and \eqref{cond4} is given by
 \begin{equation}
d = u_3v_4 + u_4v_3,\quad  c = u_3v_4 - u_4v_3, \quad e = -2u_3v_3,\quad  x_2 = -u_3v_3 - u_4v_4, \label{soleq34}
\end{equation}
where $u_3, u_4, v_3, v_4$ are arbitrary parameters.

The two solutions \eqref{soleq12} and \eqref{soleq34} will be consistent if the values of $e$ given by the two solutions are the same, that is, $u_1v_1=u_3v_3$, and hence we take 
\begin{equation}
u_1=m_1n_1, \quad v_1=m_2n_2, \quad u_3=m_1n_2, \quad v_3=m_2n_1, \label{valuv}
\end{equation}
where $m_1, m_2, n_1$ and $n_2$ are arbitrary nonzero rational parameters. On substituting the values of $u_i, v_i, i=1, 2$ given by \eqref{valuv} in \eqref{soleq12} and \eqref{soleq34}, we get the following relations,
\begin{equation}
\begin{aligned}
a & =  m_1n_1v_2 + m_2n_2u_2, \\
 b & =  m_1n_1v_2 - m_2n_2u_2, \\
c &= m_1n_2v_4 - m_2n_1u_4, \\
 d &= m_1n_2v_4 + m_2n_1u_4,\\
e & =  -2m_1m_2n_1n_2, \\
x_1 & =  -m_1m_2n_1n_2 - u_2v_2, \\
 x_2 & =  -m_1m_2n_1n_2 - u_4v_4, 
\end{aligned}
\label{valabcd}
\end{equation}
where $m_1, m_2, n_1, n_2, u_2, u_4, v_2$ and $ v_4$ are arbitrary parameters.

On substituting the values of $a$ and $x_1$ given by \eqref{valabcd} in \eqref{cond1}, we get $y_1^2=-(m_2n_2 - v_2)(m_2n_2 + v_2)(m_1n_1 - u_2)(m_1n_1 + u_2)$. Thus the right-hand side of this last equation must be a nonzero perfect square, and this is possible if and only if there exists a nonzero rational number $h_1$ such that $h_1^2(m_1n_1 - u_2)(m_1n_1 + u_2)=-(m_2n_2 - v_2)(m_2n_2 + v_2)$, hence  there exist nonzero rational parameters $p_i, q_i, i=1, 2$, such that
\begin{equation}
\begin{aligned}
h_1(m_1n_1 - u_2)& = 2p_1q_2, \quad &h_1(m_1n_1 + u_2) & = 2p_2q_1, \\
 -(m_2n_2 - v_2)& = 2p_1q_1, \quad & m_2n_2 + v_2& = 2p_2q_2. 
\end{aligned}
\label{condmn12}
\end{equation}
On solving Eqs. \eqref{condmn12} for $m_1, m_2, u_2$ and $v_2$, we get,
\begin{equation}
\begin{aligned}
m_ 1 & =  (p_ 1q_ 2 + p_ 2q_ 1)/(h_1n_ 1), & m_ 2 & =  -(p_ 1q_ 1 - p_ 2q_ 2)/n_ 2, \\
u_2 & =  -(p_ 1q_ 2 - p_ 2q_ 1)/h_1, & v_2 & =  p_ 1q_ 1 + p_ 2q_ 2.
\end{aligned}
\label{valmn12}
\end{equation}

Similarly, on substituting the values of $d$ and $x_2$ given by \eqref{valabcd} in \eqref{cond4}, we get $y_2^2=-(m_2n_1 - v_4)(m_2n_1 + v_4)(m_1n_2 - u_4)(m_1n_2 + u_4)$, and proceeding as in the previous paragraph, we get the solution,
\begin{equation}
\begin{aligned}
m_1 &= (r_1s_2 + r_2s_1)/(h_2n_2), & m_2 &= -(r_1s_1 - r_2s_2)/n_1, \\
u_4 &= -(r_1s_2 - r_2s_1)/h_2, & v_4 &= r_1s_1 + r_2s_2,
\end{aligned}
\label{valmn12sec}
\end{equation}
where $h_2, r_1, r_2, s_1, s_2$ are arbitrary nonzero rational parameters.

For the solutions \eqref{valmn12} and \eqref{valmn12sec} to be consistent, we must impose the conditions that the respective values of $m_1$ and $m_2$ given by these two solutions coincide. These conditions are satisfied if and only if we choose $h_1, h_2, n_1$ and $n_2$ as follows:
\begin{equation}
\begin{aligned}
 h_1 & =  (p_1q_2 + p_2q_1)(p_1q_1 - p_2q_2), & h_2 & =  (r_1s_2 + r_2s_1)(r_1s_1 - r_2s_2),\\
 n_1 & =  r_1s_1 - r_2s_2, & n_2 & =  p_1q_1 - p_2q_2.
\end{aligned}
\label{valhn}
\end{equation}

Using the values of $h_i, n_i, i=1, 2$, given by \eqref{valhn} we can now find  the values of $m_1, m_2, u_2, v_2, u_4$ and $ v_4$ from the relations \eqref{valmn12} and \eqref{valmn12sec}, and we can  solve Eqs. \eqref{cond1} and \eqref{cond4} to obtain rational values of $y_1, y_2$, and further,  using the relations \eqref{valabcd}, we get, on appropriate scaling, the following solution of the simultaneous diophantine Eqs. \eqref{cond1}--\eqref{cond4}: 
\begin{equation}
\begin{aligned}
 a & =  q_1q_2(r_1s_2 + r_2s_1)(r_1s_1 - r_2s_2)(p_1^2 + p_2^2),\\
 b & =  p_1p_2(r_1s_2 + r_2s_1)(r_1s_1 - r_2s_2)(q_1^2 + q_2^2),\\
 c & =  r_1r_2(p_1q_2 + p_2q_1)(p_1q_1 - p_2q_2)(s_1^2 + s_2^2),\\
 d & =  s_1s_2(p_1q_2 + p_2q_1)(p_1q_1 - p_2q_2)(r_1^2 + r_2^2),\\
 e & =  (p_1q_2 + p_2q_1)(p_1q_1 - p_2q_2)(r_1s_2 + r_2s_1)(r_1s_1 - r_2s_2),\\
 x_1 & = q_1q_2(r_1s_2 + r_2s_1)(r_1s_1 - r_2s_2)(p_1 - p_2)(p_1 + p_2), \\
x_2 & = s_1s_2(p_1q_2 + p_2q_1)(p_1q_1 - p_2q_2)(r_1 - r_2)(r_1 + r_2),\\
y_1 & =  2p_ 1p_ 2q_1q_ 2(r_ 1s_ 2 + r_ 2s_ 1)(r_ 1s_ 1 - r_ 2s_ 2), \\
 y_2 & =  -2r_ 1r_ 2s_1s_ 2(p_ 1q_ 2 + p_ 2q_ 1)(p_ 1q_ 1 - p_ 2q_ 2).
\end{aligned}
\label{valabcdequadone}
\end{equation}
where $p_i, q_i, r_i,  s_i, i=1, 2$, are arbitrary rational parameters. 

While solving the simultaneous diophantine equations \eqref{cond1}--\eqref{cond4}, we had to solve several intermediate equations and for each of these equations, we obtained the complete solution. Hence, the formulae  \eqref{valabcdequadone} give the complete solution of the simultaneous  equations \eqref{cond1}--\eqref{cond4}.

It follows that the formulae \eqref{valabcdequadone}  generate all convex quadrilatetrals whose four sides $a, b, c, d$, one diagonal and the area are given by rational numbers. We must, of course, choose the parameters such that the rational numbers mentioned towards the end of Theorem \ref{Thcond} are all positive. The area of the quadrilateral defined by \eqref{valabcdequadone} is given by  
\begin{multline}
A= (p_1q_2 + p_2q_1)(p_1q_1 - p_2q_2)(r_1s_2 + r_2s_1)(r_1s_1 - r_2s_2)\\
\times [p_1^2q_1q_2r_1r_2s_1s_2 + \{q_1^2r_1r_2s_1s_2 + (r_1s_2 + r_2s_1)(r_1s_1 - r_2s_2)q_1q_2\\
 - q_2^2r_1r_2s_1s_2\}p_1p_2 - p_2^2q_1q_2r_1r_2s_1s_2] . \label{areaquadone}
\end{multline}

\subsection{Rational quadrilaterals}\label{rationalquad}
We will now obtain parametrizations of rational quadrilaterals by obtaining parametric solutions of the simultaneous diophantine Eqs. \eqref{cond1}--\eqref{cond5}. 

On substituting the values  of $x_1, x_2, y_1$ and $y_2$ given by \eqref{valabcdequadone} in Eq. \eqref{cond5}, we get the condition, 
\begin{multline}
q_1^2q_2^2r_1^2r_2^2(p_1^2 + p_2^2)^2s_1^4 - 2p_1p_2q_1q_2r_1r_2\{(q_1r_1 + q_1r_2 + q_2r_1 - q_2r_2)p_1 + (q_1r_1 - q_1r_2\\
 - q_2r_1 - q_2r_2)p_2\}\{(q_1r_1 - q_1r_2 - q_2r_1 - q_2r_2)p_1 - (q_1r_1 + q_1r_2 + q_2r_1 - q_2r_2)p_2\}s_1^3s_2\\
 + \{2p_1^4q_1^2q_2^2r_1^2r_2^2 + 8q_1q_2r_1r_2(q_1r_2 + q_2r_1)(q_1r_1 - q_2r_2)p_1^3p_2 + (q_1^4r_1^4 + 2q_1^4r_1^2r_2^2 + q_1^4r_2^4\\
 + 8q_1^3q_2r_1^3r_2 - 8q_1^3q_2r_1r_2^3 + 2q_1^2q_2^2r_1^4 - 24q_1^2q_2^2r_1^2r_2^2 + 2q_1^2q_2^2r_2^4 - 8q_1q_2^3r_1^3r_2 \\
+ 8q_1q_2^3r_1r_2^3 + q_2^4r_1^4 + 2q_2^4r_1^2r_2^2 + q_2^4r_2^4)p_1^2p_2^2 - 8q_1q_2r_1r_2(q_1r_2 + q_2r_1)(q_1r_1 - q_2r_2)p_1p_2^3\\
 + 2p_2^4q_1^2q_2^2r_1^2r_2^2\}s_1^2s_2^2 + 2p_1p_2q_1q_2r_1r_2\{(q_1r_1 + q_1r_2 + q_2r_1 - q_2r_2)p_1 + (q_1r_1 - q_1r_2 \\
- q_2r_1 - q_2r_2)p_2\}\{(q_1r_1 - q_1r_2 - q_2r_1 - q_2r_2)p_1 - (q_1r_1 + q_1r_2 + q_2r_1 - q_2r_2)p_2\}s_1s_2^3\\
 + q_1^2q_2^2r_1^2r_2^2(p_1^2 + p_2^2)^2s_2^4=f^2. \label{cond5f}
\end{multline}

We need to solve Eq. \eqref{cond5f} such that the  rational numbers mentioned towards the end of Theorem \ref{Thcond} are positive. Thus the parameters $p_i, q_i, r_i, s_i, i=1, 2$, must all be nonzero and simple solutions of Eq. \eqref{cond5f} such as $p_1=p_2q_2/q_1$ or $p_1=-p_2q_1/q_2$ or $r_1=r_2s_2/s_1$ or $r_1=-r_2s_1/s_2$ are also ruled out. We will give four simple solutions of \eqref{cond5f} that satisfy the conditions of Theorem \ref{Thcond}, and thus obtain four parametrizations of rational quadrilaterals. We will also show how infinitely many parametrizations of rational quadrilaterals may be obtained.

\subsubsection{Rational cyclic quadrilaterals}
We first give a solution of Eq. \eqref{cond5f} that generates all rational cyclic quadrilaterals. 
It is readily verified that if we take
\begin{equation}
\begin{aligned}
s_1&=p_1(q_1r_2 + q_2r_1) + p_2(q_1r_1 - q_2r_2),\\
 s_2&=p_1(q_1r_1 - q_2r_2) - p_2(q_1r_2 + q_2r_1),
\end{aligned}
\label{sol1cond5}
\end{equation}
the left-hand side of Eq. \eqref{cond5f} becomes a perfect square. This yields a rational quadrilateral whose    sides $a, b, c, d$, the  diagonals $e, f$ and the area $A$ may be written, after appropriate scaling, as follows: 
\begin{equation}
\begin{aligned}
a  & =  (p_1^2 + p_2^2)(r_1^2 + r_2^2)q_1q_2,\\
 b & = (q_1^2 + q_2^2)(r_1^2 + r_2^2)p_1p_2, \\
 c & = (q_1^2 + q_2^2)(p_1^2 + p_2^2)r_1r_2,\\
 d  & =  (p_1q_1r_2 + p_1q_2r_1 + p_2q_1r_1 - p_2q_2r_2)\\
& \quad \quad  \times (p_1q_1r_1 - p_1q_2r_2 - p_2q_1r_2 - p_2q_2r_1),\\
 e  & = (r_1^2 + r_2^2)(p_1q_2 + p_2q_1)(p_1q_1 - p_2q_2), \\
f  & =  (q_1^2 + q_2^2)(p_1r_2 + p_2r_1)(p_1r_1 - p_2r_2), \\
A & = (p_1q_1 - p_2q_2)(p_1q_2 + p_2q_1)(p_1r_1 - p_2r_2)\\
& \quad \quad  \times (p_1r_2 + p_2r_1)(q_1r_1 - q_2r_2)(q_1r_2 + q_2r_1).
\end{aligned}
\label{parmsol1}
\end{equation}

Further, the related values of $y_1, y_2$, obtained from \eqref{valabcdequadone} are as follows:
\begin{equation}
\begin{aligned}
x_1& = (p_1 + p_2)(p_1 - p_2)(r_1^2 + r_2^2)q_1q_2,\\
x_2&=(r_1 - r_2)(r_1 + r_2)(p_1q_1r_2 + p_1q_2r_1 + p_2q_1r_1 - p_2q_2r_2)\\
& \quad \quad  \times (p_1q_1r_1 - p_1q_2r_2 - p_2q_1r_2 - p_2q_2r_1)/(r_1^2 + r_2^2),\\
y_1 & =  2(r_1^2 + r_2^2)p_1p_2q_1q_2, \\
y_2 & =  -2r_1r_2(p_1q_1r_2 + p_1q_2r_1 + p_2q_1r_1 - p_2q_2r_2)\\
& \quad \quad  \times (p_1q_1r_1 - p_1q_2r_2 - p_2q_1r_2 - p_2q_2r_1)/(r_1^2 + r_2^2),
\end{aligned}
\label{valyparmsol1}
\end{equation}

If  we assign positive rational values to all  the parameters $p_i, q_i, r_i, i=1, 2$, such that $q_1 > q_2r_2/r_1$ and $p_1 > p_2(q_1r_2 + q_2r_1)/(q_1r_1 - q_2r_2)$, it readily follows that  all  conditions mentioned in Theorem \ref{Thcond} are satisfied so that the formulae \eqref{parmsol1} define a rational convex quadrilateral. As a numerical example,  when we take  $(p_1, p_2, q_1, q_2, r_1, r_2)=(4, 1, 3, 1, 2, 1)$, we get a rational convex quadrilateral whose sides, diagonals and the area, after appropriate scaling, are given by $(a, b, c, d, e, f, A) = (51, 40, 68, 75, 77, 84, 3234)$.

We will now show that  the formulae \eqref{parmsol1} give the sides, diagonals and the area of a cyclic quadrilateral. In fact, it is useful to recall here that if $a, b, c, d$ are the consecutive sides of a cyclic quadrilateral, the area $K$ and  the lengths $e$ and $f$   of the two diagonals  of the quadrilateral are given by the following formulae proved by Brahmagupta \cite[p. 187]{Ev}:
\begin{equation}
K=\sqrt{(s-a)(s-b)(s-c)(s-d)}, \label{Brahmaarea}
\end{equation}
and 
\begin{equation}
\begin{aligned}
e&=\sqrt{(ab+cd)(ac+bd)/(ad+bc)},\\
f&=\sqrt{(ac+bd)(ad+bc)/(ab+cd)}, 
\end{aligned}
\label{Brahmadiag}
\end{equation}
where $s$ is the semi-perimeter, that is,
\begin{equation}
s=(a+b+c+d)/2. \label{vals}
\end{equation}
Further,  the circumradius $R$  of the cyclic quadrilateral is given by the following formula of Paramesvara \cite{Gu}:
\begin{equation}
R=\frac{\displaystyle 1}{\displaystyle 4}\ \sqrt{\frac{\displaystyle(ab+cd)(ac+bd)(ad+bc)}{\displaystyle (s-a)(s-b)(s-c)(s-d)}} .
\end{equation}

We also note that the area $\Delta$ of an arbitrary convex quadrilateral having sides $a, b, c, d$, and the sum of any one pair of opposite angles equal to $2u$,  is given by the following formula \cite[p. 82]{Jo}:
\begin{equation}
\Delta^2 =(s-a)(s-b)(s-c)(s-d) - abcd \cos^2 u, \label{areagenquad}
\end{equation}
where $s$ is the semiperimeter.

It follows from the formulae \eqref{Brahmaarea} and \eqref{areagenquad} that a quadrilateral with sides $a, b, c, d$, is cyclic if and only if its area is given by the formula \eqref{Brahmaarea}. It is readily verified that the area $A$ of the quadrilateral defined by \eqref{parmsol1} is exactly the area of a cyclic quadrilateral whose sides $a, b, c, d$, are given by \eqref{parmsol1}. Further, the lengths of the two diagonals, given by the formulae \eqref{Brahmadiag}, coincide with the lengths $e$ and $f$ of the two diagonals given by \eqref{parmsol1}. Thus, the quadrilateral defined by \eqref{parmsol1} is a cyclic quadrilateral.  We also note that that the circumradius of the quadrilateral defined by \eqref{parmsol1} is
\begin{equation}
R=(p_1^2 + p_2^2)(q_1^2 + q_2^2)(r_1^2 + r_2^2)/4.
\end{equation}

It is pertinent to consider at this stage whether there are any other values of the parameters $p_i, q_i, r_i, s_i, i=1, 2$, such the quadrilateral defined by \eqref{valabcdequadone} and \eqref{areaquadone} becomes a cyclic quadrilateral. The condition that the area $A$ of this quadrilateral becomes equal to the area of a cyclic quadrilateral with sides $a, b, c, d$, may be written as follows:
\begin{multline}
p_1p_2q_1q_2r_1r_2s_1s_2(r_1s_2 + r_2s_1)^2(r_1s_1 - r_2s_2)^2(p_1q_2 + p_2q_1)^2\\
(p_1q_1 - p_2q_2)^2\{(p_1q_1r_1 - p_1q_2r_2 - p_2q_1r_2 - p_2q_2r_1)s_1\\
 -(p_1q_1r_2 + p_1q_2r_1 + p_2q_1r_1 - p_2q_2r_2)s_2\}^2=0. \label{condcyclic}
\end{multline}
Equating to $0$ any factor on the left-hand side of \eqref{condcyclic}, except the last, does not yield a quadrilateral since at least one of the four sides $a, b, c, d$, becomes $0$, and the last factor, when equated to $0$, yields the cyclic quarilateral given by \eqref{parmsol1}. 

Even if we consider the possibility that for a certain choice of parameters,   the formulae \eqref{valabcdequadone} may yield negative values for one or more of the numbers $a, b, c, d$, and we  obtain a rational convex quadrilateral by changing the signs of $a, b, c, d$, as necessary, when  we impose the condition that the area  of the quadrilateral becomes equal to the area of a cyclic quadrilateral, and proceed as above, we eventually get the same  parametrization of a quadrilateral as is given by \eqref{parmsol1}. It follows that the  formulae \eqref{parmsol1} give the complete parametrization of all cyclic quadrilaterals.
 
We note that Euler (as quoted by  Dickson \cite[p. 221]{Di2})  had given  complicated expressions for the sides and diagonals of a  cyclic quadrilateral to be rational but  the area is not  rational.  As mentioned in the Introduction, Sastry \cite[pp.\ 170--171]{Sa} has given a parametrization of cyclic quadrilaterals but the formulae \eqref{parmsol1} are simpler and more symmetric as compared to Sastry's formulae.

\subsubsection{A parametrization of noncyclic convex quadrilaterals} 
A second solution of Eq. \eqref{cond5f} is obtained if we take
\[s_1 = r_1r_2\{2p_1q_1q_2 + (q_1^2 - q_2^2)p_2\}, \quad  s_2 =q_1q_2\{(r_1^2 - r_2^2)p_1 - 2p_2r_1r_2\}.\]
This yields a rational quadrilateral whose sides $a, b, c, d$, the  diagonals $e, f$ and the area $A$ are as follows:
\begin{equation}
\begin{aligned}
a& =q_1q_2(p_1^2 + p_2^2)\{(r_1^2 + r_2^2)p_1q_1q_2 + (q_1^2r_2 - 2q_1q_2r_1 - q_2^2r_2)p_2r_2\}\\
& \quad \quad \times \{(r_1^2 + r_2^2)p_1q_1q_2 + (q_1^2r_1 + 2q_1q_2r_2 - q_2^2r_1)p_2r_1\},\\
b& =p_1p_2(q_1^2 + q_2^2)\{(r_1^2 + r_2^2)p_1q_1q_2 +(q_1^2r_2 - 2q_1q_2r_1 - q_2^2r_2)p_2 r_2\}\\
& \quad \quad \times \{(r_1^2 + r_2^2)p_1q_1q_2 + (q_1^2r_1 + 2q_1q_2r_2 - q_2^2r_1)p_2r_1\},\\
c& =(p_1q_2 + p_2q_1)(p_1q_1 - p_2q_2)\{(r_1^2 + r_2^2)^2p_1^2q_1^2q_2^2 + 4(q_1r_1+ q_2r_2)q_1q_2r_1r_2\\
& \quad \quad \times (q_1r_2 - q_2r_1)p_1p_2 + (q_1^2 + q_2^2)^2p_2^2r_1^2r_2^2\},\\
d& =q_1q_2(p_1q_2 + p_2q_1)(p_1q_1 - p_2q_2)(r_1^2 + r_2^2)\{2p_1q_1q_2 + (q_1^2 - q_2^2)p_2\}\\
& \quad \quad \times \{(r_1^2 - r_2^2)p_1 - 2p_2r_1r_2\}, \\
e& =(p_1q_2 + p_2q_1)(p_1q_1 - p_2q_2)\{(r_1^2 + r_2^2)p_1q_1q_2 + (q_1^2r_2 - 2q_1q_2r_1\\ 
& \quad \quad  - q_2^2r_2)p_2r_2\}\{(r_1^2 + r_2^2)p_1q_1q_2 + (q_1^2r_1 + 2q_1q_2r_2 - q_2^2r_1)p_2r_1\},\\
f& =q_1q_2\{(r_1^2 + r_2^2)^2p_1^4q_1^2q_2^2 + 2(q_1^2r_1^4 + q_1^2r_2^4 + 2q_1q_2r_1^3r_2 \\
&\quad \quad - 2q_1q_2r_1r_2^3 - q_2^2r_1^4 - q_2^2r_2^4)p_1^3p_2q_1q_2 + (q_1^4r_1^4 - q_1^4r_1^2r_2^2 + q_1^4r_2^4\\
&\quad \quad - q_1^2q_2^2r_1^4 - 8q_1^2q_2^2r_1^2r_2^2 - q_1^2q_2^2r_2^4 + q_2^4r_1^4 - q_2^4r_1^2r_2^2 + q_2^4r_2^4)p_1^2p_2^2\\
&\quad \quad - 2(q_1^4r_1^2 - q_1^4r_2^2 + 2q_1^3q_2r_1r_2 - 2q_1q_2^3r_1r_2 + q_2^4r_1^2 - q_2^4r_2^2)p_1p_2^3r_1r_2\\
 &\quad \quad + (q_1^2 + q_2^2)^2p_2^4r_1^2r_2^2\},\\
A& =q_1q_2(p_1q_2 + p_2q_1)(p_1q_1 - p_2q_2)\{(r_1^2 + r_2^2)p_1q_1q_2 +(q_1^2r_2 - 2q_1q_2r_1 \\
& \quad \quad - q_2^2r_2)p_2 r_2\}\{(r_1^2 + r_2^2)p_1q_1q_2 + (q_1^2r_1 + 2q_1q_2r_2 - q_2^2r_1)p_2r_1\}\\
& \quad \quad \times \{2p_1^2q_1q_2r_1r_2 + (q_1r_2 + q_2r_1)(q_1r_1 - q_2r_2)p_1p_2 - 2p_2^2q_1q_2r_1r_2\}\\
& \quad \quad \times \{(r_1^2 - r_2^2)p_1^2q_1q_2 + (r_1^2 - r_2^2)(q_1^2 - q_2^2)p_1p_2 - (q_1^2 - q_2^2)p_2^2r_1r_2\},
\end{aligned}
\label{parmsol2}
\end{equation}
where $p_i, q_i, r_i, i=1, 2$, are arbitrary nonzero rational parameters.

As a numerical example, if we take  $(p_1, p_2, q_1, q_2, r_1, r_2)=(3, 1, 2, 1$, $3, 1)$, we get a rational convex quadrilateral whose sides, diagonals and the area, after appropriate scaling, are given by
$(a, b, c, d, e, f, A) = (748, 561, 615$, $1000, 935, 1068, 490314)$.

We note that if, in the formulae \eqref{parmsol2}, we take $q_1, q_2, r_1, r_2$ as numerical constants, the sides and diagonals of the quadrilateral are given by quartic polynomials in the parameters $p_1$ and $p_2$. Similarly, we may consider the formulae for the sides and diagonals of the above quadrilateral as  quartic polynomials in the parameters $r_1$ and $r_2$.

\subsubsection{A second parametrization of noncyclic convex quadrilaterals} 
To obtain a second  parametrization of rational noncyclic convex quadrilaterals, we note that the left-hand side of Eq. \eqref{cond5f} becomes a perfect square if we take 
\[
r_1=-p_2,\; r_2=p_1, \;s_1=(p_1-p_2)q_1-(p_1+p_2)q_2, \; s_2=-(p_1+p_2)q_1-(p_1-p_2)q_2. \]
We thus get a rational quadrilateral whose sides $a, b, c, d$, the  diagonals $e, f$ and the area $A$ may be written, after replacing the parameters $p_2$ and $q_2$ by $-p_2$ and $-q_2$ respectively,  as follows: 
\begin{equation}
\begin{aligned}
a & = (p_1^2 + p_2^2)^2(q_2^2 - q_1^2)q_1q_2, \\
b & =  p_1p_2(p_1^2 + p_2^2)(q_2^4  - q_1^4), \\
c & =  2p_1p_2(p_1q_2 + p_2q_1)(p_2q_2-p_1q_1)(q_1^2 + q_2^2),\\
 d & =  (p_1q_2 + p_2q_1)(p_2q_2 - p_1q_1)(p_1q_1 - p_1q_2 - p_2q_1 - p_2q_2)\\
& \quad \quad \times (p_1q_1 + p_1q_2 + p_2q_1 - p_2q_2), \\
e & =  (p_1^2 + p_2^2)(p_2q_2 - p_1q_1)(p_1q_2 + p_2q_1)(q_2^2- q_1^2), \\
f& = p_1p_2(p_2^2-p_1^2)(q_1^2 + q_2^2)^2,\\
A& = p_1p_2(p_2^2 - p_1^2)(q_2^2 - q_1^2)(p_1q_2 + p_2q_1)(p_2q_2 - p_1q_1)\\
& \quad \quad \times (p_1q_1^2 - p_1q_2^2 - 2p_2q_1q_2)(2p_1q_1q_2 + p_2q_1^2 - p_2q_2^2),
\end{aligned}
\label{parmsol3}
\end{equation}
where $p_i, q_i, i=1, 2$, are arbitrary nonzero rational parameters.

The values of $a, b, c, d, e$ and $f$ given by \eqref{parmsol3} and the values of $x_1, x_2, y_1$, $y_2$ given by 
 \begin{equation}
\begin{aligned}
x_1 & = (p_2^4  - p_1^2)(q_2^2- q_1^2)q_1q_2,\\
 x_2 & = (p_1q_1 - p_1q_2 - p_2q_1 - p_2q_2)(p_1q_1 + p_1q_2 + p_2q_1 - p_2q_2)(p_2^2 - p_1^2)\\
& \quad \quad \times (p_1q_2 + p_2q_1)(p_2q_2 - p_1q_1)/(p_1^2 + p_2^2),\\
 y_1 & = 2(p_1^2 + p_2^2)(q_2^2 - q_1^2)p_1p_2q_1q_2,\\
 y_2 & = -2p_1p_2(p_1q_1 - p_1q_2 - p_2q_1 - p_2q_2)(p_1q_1 + p_1q_2 + p_2q_1 - p_2q_2)\\
& \quad \quad \times (p_1q_2 + p_2q_1)(p_2q_2 - p_1q_1)/(p_1^2 + p_2^2)
\end{aligned}
\label{parmsol3xy}
\end{equation}
satisfy the simultaneous Eqs. \eqref{cond1}--\eqref{cond5}.

It follows from Theorem \ref{Thcond} that if we assign positive rational values to the parameters $p_i, q_i, i=1, 2$, such that $q_2 > q_1$, and $p_2 > p_1(q_1 + q_2)/(q_2 - q_1)$, the quadrilateral defined by \eqref{parmsol3} is a convex quadrilateral.

As a numerical example, on taking  $(p_1, p_2, q_1, q_2)=(1, 2, 1, 5)$,  we get a rational convex quadrilateral whose sides, diagonals and the area, after appropriate scaling, are given by
$(a, b, c, d, e, f, A) =(125, 260, 273, 84, 315, 169$, $26334)$.

\subsubsection{Rational quadrilaterals with two equal sides}
There are special quadrilaterals such as paralellograms, rhombi amd kites in which two of the four sides are equal in addition to  other geometrical conditions. Parametrizations of such special rational quadrilaterals may be obtained by considering the problem ab initio rather than by considering it as a special case of the general problem of finding rational quadrilaterals.  For instance, the problem of finding a rational paralellogram is equivalent to finding a rational  triangle with a rational median, and its complete solution has been given by Dickson \cite[pp. 249--250]{Di1}. We will not be concerned with such special quadrilaterals in this paper.  We  give below a parametrization  of a more general quadrilateral in which two of the four sides of the quadrilateral are necessarily equal.

We note that a solution of Eq. \eqref{cond5f} may be obtained by taking
\[q_1=-p_2,\; q_2=p_1,\; s_1=(p_1-p_2)r_1-(p_1+p_2)r_2,\; s_2=-(p_1+p_2)r_1-(p_1-p_2)r_2. \]
This yields  a rational quadrilateral whose sides $a, b, c, d$, the  diagonals $e, f$ and the area $A$ are given by the following formulae:
\begin{equation}
\begin{aligned}
a & = b  = (p_1^2 + p_2^2)(r_1^2 + r_2^2), \quad  c  =  4(p_1^2 + p_2^2)r_1r_2, \\
d & =  2(p_1r_1 + p_1r_2 + p_2r_1 - p_2r_2)(p_1r_1 - p_1r_2 - p_2r_1 - p_2r_2),  \\
 e & =  2(r_1^2 + r_2^2)(p_2^2-p_1^2),  \quad f =  (p_1^2 + p_2^2)(r_1^2 + r_2^2),\\
A& = 2(p_2^2 - p_1^2)(2p_1r_1r_2 + p_2r_1^2 - p_2r_2^2)(p_1r_1^2 - p_1r_2^2 - 2p_2r_1r_2),
\end{aligned}
\label{parmsol4}
\end{equation}
where $p_i, r_i, i=1, 2$, are arbitrary parameters.

It follows from Theorem \ref{Thcond} that if we assign positive rational values to the parameters $p_i, q_i, i=1, 2$, such that $r_2 > r_1$ and $p_2 > p_1(r_1+r_2)/(r_2-r_1)$, the formulae \eqref{parmsol4} will always generate rational convex quadrilaterals with two equal sides.
 As a numerical example, on taking $(p_1, p_2, r_1, r_2)= (1, 3, 1, 3)$, we get  a rational convex quadrilateral whose sides, diagonals and the area, after appropriate scaling, are given by
$(a, b, c, d, e, f, A) =(25, 25, 30, 14, 40, 25, 468)$.

\subsubsection{Further parametrizations of rational quadrilaterals}
We will now show how additional parametrizations of rational quadrilaterals may be obtained by finding more solutions of the diophantine Eq. \eqref{cond5f}. On writing $s_1=s_2X, f=s_2^2Y$, Eq. \eqref{cond5f} reduces to
\begin{multline}
Y^2=q_1^2q_2^2r_1^2r_2^2(p_1^2 + p_2^2)^2X^4 - 2p_1p_2q_1q_2r_1r_2\{(q_1r_1 + q_1r_2 + q_2r_1 - q_2r_2)p_1\\
 + (q_1r_1 - q_1r_2 - q_2r_1 - q_2r_2)p_2\}\{(q_1r_1 - q_1r_2 - q_2r_1 - q_2r_2)p_1 - (q_1r_1 + q_1r_2\\
 + q_2r_1 - q_2r_2)p_2\}X^3+ \{2p_1^4q_1^2q_2^2r_1^2r_2^2 + 8q_1q_2r_1r_2(q_1r_2 + q_2r_1)(q_1r_1 - q_2r_2)p_1^3p_2\\
 + (q_1^4r_1^4 + 2q_1^4r_1^2r_2^2 + q_1^4r_2^4+ 8q_1^3q_2r_1^3r_2 - 8q_1^3q_2r_1r_2^3 + 2q_1^2q_2^2r_1^4 - 24q_1^2q_2^2r_1^2r_2^2\\ + 2q_1^2q_2^2r_2^4 - 8q_1q_2^3r_1^3r_2 + 8q_1q_2^3r_1r_2^3 + q_2^4r_1^4 + 2q_2^4r_1^2r_2^2 + q_2^4r_2^4)p_1^2p_2^2 \\
- 8q_1q_2r_1r_2(q_1r_2 + q_2r_1)(q_1r_1 - q_2r_2)p_1p_2^3+ 2p_2^4q_1^2q_2^2r_1^2r_2^2\}X^2 + 2p_1p_2q_1q_2r_1r_2\\
 \times \{(q_1r_1 + q_1r_2 + q_2r_1 - q_2r_2)p_1 + (q_1r_1 - q_1r_2 - q_2r_1 - q_2r_2)p_2\}\{(q_1r_1 - q_1r_2 - q_2r_1\\ - q_2r_2)p_1
 - (q_1r_1 + q_1r_2 + q_2r_1 - q_2r_2)p_2\}X+ q_1^2q_2^2r_1^2r_2^2(p_1^2 + p_2^2)^2. \label{cond5XY}
\end{multline}

The right-hand side of Eq. \eqref{cond5XY} is a quartic polynomial in $X$, and Eq. \eqref{cond5XY} may be considered as representing the quartic model of an elliptic curve over the function field $\mathbb{Q}(p_1. p_2, q_1, q_2, r_1, r_2)$. A point $P=(X_1, Y_1)$ on the elliptic curve \eqref{cond5XY}, derived from the solution \eqref{sol1cond5} of Eq. \eqref{cond5f}, is as follows:
\[
\begin{aligned}
(X_1, Y_1) & =(\{(q_1r_2 + q_2r_1)p_1 + (q_1r_1 - q_2r_2)p_2\}\\
& \quad \quad \times \{(q_1r_1 - q_2r_2)p_1 + (-q_1r_2 - q_2r_1)p_2\}^{-1},\\
& \quad \quad (p_1q_2 + p_2q_1)(p_1q_1 - p_2q_2)(p_1r_2 + p_2r_1)(p_1r_1 - p_2r_2)(q_1^2 + q_2^2)\\
& \quad \quad \times (r_1^2 + r_2^2)\{(q_1r_1 - q_2r_2)p_1 + (-q_1r_2 - q_2r_1)p_2\}^{-2}).
\end{aligned}
\]

Since the coefficient of $X^4$ in the quartic polynomial on the right-hand side of \eqref{cond5XY} is a perfect square, we can, by following a well-known procedure (see, for instance, \cite[p. 77]{Mo2} or \cite[pp. 35--36]{Ca}),  transform the elliptic curve \eqref{cond5XY} to the usual cubic model, and find a point $P_1$ on the cubic model  corresponding to the point $P$ on the quartic model given by \eqref{cond5XY}. It is easily established  that the point $P_1$ is not of finite order. Thus, by repeated application of the group law, we can find infinitely many rational points on the cubic model of the curve, and hence also on the curve \eqref{cond5XY}. We thus get infinitely many rational solutions of Eq. \eqref{cond5XY}, and we can thus find infinitely  many parametrizations of rational quadrilaterals.

\section{Concluding remarks}
In this paper we described a method  of generating all  quadrilaterals whose sides, diagonals and the area are given by rational numbers. This method is different from the classic construction of rational quadrilaterals given by Kummer in 1848. We were able to obtain a complete parametrization of all rational cyclic quadrilaterals.  For general rational quadrilaterals, we obtained several parametrizations in terms of quartic polynomials. When two  sides of the quadrilateral are equal, we obtained  a parametrization in terms of quadratic polynomials. The formulae for  rational quadrilaterals obtained in this paper are simpler than the known parametrizations of such quadrilaterals.

\noindent Postal address: Ajai Choudhry, 13/4 A, Clay Square, Lucknow - 226001, India\\

\noindent E-mail: ajaic203@yahoo.com

\end{document}